\documentclass[journal]{IEEEtran}
\usepackage{cite}
\usepackage{url}
\usepackage{empheq}
\usepackage{float}
\usepackage{dsfont}
\usepackage[hidelinks]{hyperref}

\usepackage{algorithm}
\usepackage{algorithmic}

\usepackage{amsmath,amssymb,amsfonts}
\allowdisplaybreaks[4]
\usepackage{graphicx}
\usepackage{textcomp}
\usepackage{amsmath}
\usepackage{enumerate}
\usepackage{epstopdf}
\usepackage{array}
\usepackage{booktabs}
\usepackage{subfigure}
\usepackage{multirow}
\usepackage{soul, color}
\usepackage[usenames,dvipsnames]{xcolor}
\usepackage[latin1]{inputenc}
\usepackage{cases}

\newtheorem{proposition}{Proposition}

\soulregister\cite7 
\soulregister\citep7 
\soulregister\citet7 
\soulregister\ref7 
\soulregister\pageref7

\usepackage{bbding}

\usepackage{nomencl}
\makenomenclature
\usepackage{etoolbox}
\renewcommand\nomgroup[1]{%
  \item[\bfseries
  \ifstrequal{#1}{A}{Acronyms}{%
  \ifstrequal{#1}{S}{Symbols}{%
  \ifstrequal{#1}{U}{Units}{}}}%
]}

\begin{document}

\title{A Two-Stage Online Algorithm for EV Charging Station Energy Management and Carbon Trading}

\author{
Dongxiang~Yan,
Shihan~Huang,
Sen~Li,
Xiaoyi Fan,
and Yue~Chen,~\IEEEmembership{Member,~IEEE}
\thanks{D. Yan, S. Huang, and Y. Chen are with the Department of Mechanical and Automation Engineering, the Chinese University of Hong Kong, Hong Kong SAR, China (e-mail: dongxiangyan@cuhk.edu.hk, shhuang@link.cuhk.edu.hk, yuechen@mae.cuhk.edu.hk).}
\thanks{S. Li is with the Department of Civil and Environmental Engineering, The Hong Kong University of Science and Technology, Hong Kong, China. (email: cesli@ust.hk).}
\thanks{X. Fan is with the Shenzhen Jiangxing Intelligence Inc., Shenzhen, China. (email: xiaoyifan@jiangxingai.com)}
}

\markboth{}
{Shell \MakeLowercase{\textit{et al.}}: Bare Demo of IEEEtran.cls for IEEE Journals}

\maketitle

\begin{abstract}
The increasing electric vehicle (EV) adoption challenges the energy management of charging stations (CSs) due to the large number of EVs and the underlying uncertainties. Moreover, the carbon footprint of CSs is growing significantly due to the rising charging power demand. This makes it important for CSs to properly manage their energy usage and ensure their carbon footprint stay within their carbon emission quotas. This paper proposes a two-stage online algorithm for this purpose, considering the different time scales of energy management and carbon trading. In the first stage, the CS characterizes the real-time aggregate EV power flexibility, in terms of upper and lower bounds on the total charging power, by a Lyapunov optimization-based online algorithm. In the second stage, the CS co-optimizes energy management and carbon trading, with EV charging power chosen within the aggregate flexibility region provided by the first stage. A generalized battery model is proposed to capture the dynamic carbon footprint changes and carbon trading. A virtual carbon queue is designed to develop an online algorithm for the second stage, which can ensure the carbon footprint of CS be within its carbon emission quota and its total operation cost is nearly offline optimal. Case studies validate the effectiveness and advantages of the proposed algorithm.

\end{abstract}


\IEEEpeerreviewmaketitle

\section{Introduction}
\IEEEPARstart {E}{lectric} vehicles (EVs) have achieved significant growth in recent years due to their low carbon emissions. However, the burgeoning EV population increases the complexity of energy management for charging stations (CSs). Meanwhile, though the use of EVs reduces carbon emissions in the transportation sector, considerable carbon emissions rise in the CS associated with EV charging due to the electricity consumption from the power grid \cite{dixon2020scheduling}.
To reduce operation costs and carbon footprint, it is crucial for CS operators to effectively manage the charging of EVs.

There has been extensive literature on EV charging management.
For instance, researchers have proposed dynamic EV charging strategies that track photovoltaic (PV) generation \cite{mouli2016system} and wind generation \cite{yang2018distributed} to enhance the utilization of renewable energy sources. The EV charging strategy was optimized to minimize the electricity cost given the time-of-use prices \cite{wu2017optimal}. However, the uncertainties related to EV charging were not considered. To address this issue, a probabilistic model for EV daily travel was incorporated to improve the charging strategy \cite{wu2018stochastic}.
A multi-stage energy management method was devised for a CS with solar panels and batteries \cite{yan2018optimized}.
Nevertheless, directly managing a large number of EVs is challenging since it is computational demanding.

An alternative way is to first characterize the aggregate EV power flexibility and then treat all EVs as a whole to simplify the dispatch in the subsequent stage.
Reference \cite{xu2016hierarchical} modeled the aggregate EV flexibility by directly adding the upper and lower bounds of power and energy.
References \cite{zhang2016evaluation} and \cite{wang2022evaluation} used this model to assess the vehicle-to-grid capacity of an EV fleet and to utilize EV flexibility for ensuring the reliability of power systems. Reference \cite{shi2021optimal} extended this approach by considering the spatio-temporal distribution of EVs in flexibility evaluation. In addition, an EV dispatchable region was formed to allow CS operators to participate in electricity markets \cite{zhou2021forming}.
Apart from EVs, the aggregate power flexibility of distributed energy resources in unbalanced power distribution system \cite{chen2020aggregate}, energy hubs \cite{feng2024dayahead}, and virtual power plants \cite{wang2021aggregate} has also been evaluated.

Despite the benefits of EVs, massive EV charging may raise substantial carbon footprint for CS due to the huge electricity consumption. Despite its importance, carbon emission mitigation in CS energy management has not been well explored. To reduce carbon emissions, references \cite{tulpule2013economic} and \cite{hoehne2016optimizing} analyzed the benefits of solar PV deployment and optimal charging strategy in carbon emissions reduction, respectively. A multi-objective optimization model was introduced in \cite{li2023smart} to simultaneously reduce both carbon emissions from the power grid and EV charging costs.
Cap-and-trade is an important scheme to promote carbon emission reduction \cite{yan2021blockchain}.
Charging stations whose associated carbon emissions exceed their initial free carbon emission quota must buy additional allowances from the carbon market.
The impact of cap-and-trade scheme on the operation of energy hubs \cite{zhong2023communication} and power system \cite{yang2020distributed} have been revealed, but how it influences the CS operation needs further investigation.

Moreover, the above studies adopted offline approaches that require prior knowledge of future uncertainty realizations, such as EV arrival/departure time/state-of-charge (SoC), electricity prices, and carbon intensity and prices. Obtaining such accurate information in advance can be difficult.
Therefore, online algorithms that can adapt to real-time conditions is required. Greedy algorithms were applied in real-time peer-to-peer energy market \cite{guo2021asynchronous}. However, the time-coupled constraints are neglected, resulting in suboptimal solutions. References \cite{2019zheng} and \cite{jiao2022online} utilized model predictive control methods, but they still require short-term forecasts for uncertainties.
Another technique for developing online algorithm is Lyapunov optimization. It eliminates the need for predictions and offers theoretical performance guarantees \cite{fan2020online}.
Lyapunov optimization has been applied in various domains such as microgrids \cite{2017RTMG}, shared energy storage \cite{zhong2019online}, and data centers \cite{yu2018distributed}. However, none of the above works has considered the time-coupled constraints caused by carbon footprint dynamics, which is quite different from those of EVs or energy storage. Moreover, existing works focus on online implementation with a single time scale, but EV charging and carbon trading are settled at quite different frequency. 



According to the above discussions, an online CS operation method that can adapt to uncertainties and co-optimize the energy management and carbon trading with diverse time scales is required, which is the goal of this paper. Our main contributions are two-fold:


1) \emph{Two-stage Framework}. We propose a two-stage framework for operating a CS, which can co-optimize the energy management and carbon trading while respecting their distinct time scales. In the first stage, an optimization model is developed to determine the aggregate charging power flexibility region of all EVs within the CS. We prove that any charging power trajectory within the derived region is achievable. In the second stage, the CS focuses on the energy management and carbon trading for the entire station and treats all EVs as a whole by setting the aggregate charging power within the flexibility region obtained in the first stage. A generalized battery model is proposed to capture the carbon footprint dynamics over time. This model reflects the lower frequency of carbon trading compared to energy management and is easy to extend to an online version.


2) \emph{Online Algorithm}. An online algorithm is developed to execute the two-stage framework in a prediction-free way. Specifically, we employ Lyapunov optimization to transform the offline models of the two-stage problems into their online counterparts, respectively. First, we modify the offline models by transforming the time-coupled constraints of EV charging and carbon footprint into time-average constraints. Then, charging queues and virtual carbon queues are constructed so that the time-average constraints can be further relaxed into the mean-rate-stable constraints of virtual queues. Based on these, online algorithms with feasibility and near-optimality guarantees can be derived. Compared to existing works, the proposed algorithm facilitates an online implementation with two time scales, i.e., for energy management and carbon trading, respectively.


\section{Two-stage Framework: Offline Models}
\label{secII}
\begin{figure}[htbp]
  \centering
  \includegraphics[width=0.35\textwidth]{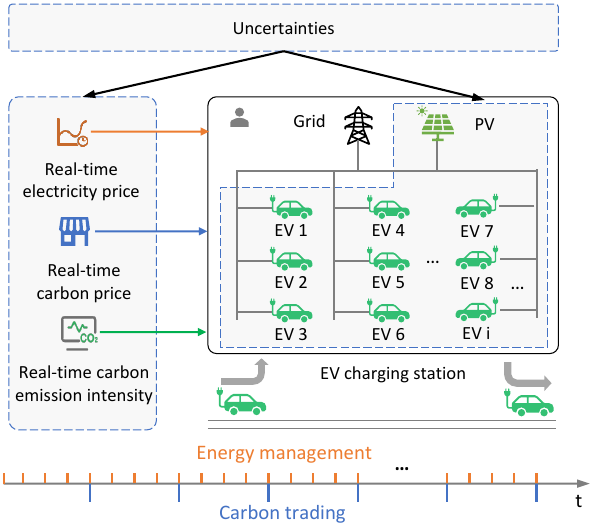}\\
  \caption{System diagram for online CS operation and carbon trading.}\label{fig:sys}
\end{figure}

Fig. \ref{fig:sys} shows the overall structure of the CS system with a large number of EVs.
Generally, the CS operator (CSO) monitors two tasks:
(1) \emph{EV charging scheduling}: When an EV, denoted by $v \in \mathcal{V}$, arrives at the CS, it submits its charging task to the CSO. Typically, $[t_v^a, t_v^d, e_v^{a}, e_v^d]$ is used to describe the charging task. Here, $t_v^a$ denotes the arrival time, $t_v^d$ the departure time, $e_v^{a}$ the initial energy level at $t_v^a$, and $e_v^d$ the intended energy level upon departure.
The CSO schedules the charging of massive EVs to ensure the charging tasks are fulfilled within the declared durations.
(2) \emph{CS energy management and carbon trading}: The CSO uses on-site renewable generation and purchases electricity from the grid to meet the charging demand at the minimum operation cost. In addition, if the carbon footprint associated with EV charging exceeds the initially allocated carbon emission quota of the CS, the CSO has to buy from the carbon market to comply with the carbon emission quota restriction. The carbon trading usually has a lower frequency than energy management.

Here, we propose a two-stage framework: In the first stage, the CSO derives the aggregate flexibility region of all EVs (i.e., feasible range of aggregate EV charging power) for every time slots; In the second stage, the CSO optimizes the real-time energy management within the charging station with the aggregate flexibility region serving as a constraint. The CSO also decides on the carbon trading with the carbon market. Energy management and carbon trading differ on the time scale as shown in Fig. \ref{fig:sys}, e.g., 10 min for energy management and 1 hour for carbon trading. A generalized battery model is proposed to characterize the time-varying carbon footprint and ensure compliance with the carbon emission quota restriction. In the following, we formulate the offline problems of the two stages, respectively. Then in Sections \ref{sec:online-evagg} and \ref{sec:online-evcs}, their online counterparts are provided.

\subsection{Stage 1: EV Aggregate Power Flexibility Characterization}
We suggest approximating the real aggregate power flexibility region of all EVs inside a CS using a sequence of intervals $[\check{p}_{d,t},\hat{p}_{d,t}],\forall t$, for convenience of use.
In other words, $[\check{p}_{d,1},\hat{p}_{d,1}] \times ... \times [\check{p}_{d,T},\hat{p}_{d,T}]$ approximates the aggregate power flexibility region and can be specified by an upper power trajectory $\{\hat p_{d,t},\forall t\}$ and a lower power trajectory $\{\check p_{d,t},\forall t\}$.
We formulate the following optimization problem: 
\begin{subequations}
\label{eq:flexibility}
\begin{align}
\textbf{P1:}~\max \limits_{\hat{p}_{d,t},\check{p}_{d,t}} ~&\mathcal{F}:=\lim\limits_{T\rightarrow\infty}\frac{1}{T}\sum\limits_{t=1}^{T}\mathbb{E}\Big[F_{t}\Big],\label{eq:p1-obj} \\
\hbox{s.t.}~ & \{\hat{p}_{d,t},\check{p}_{d,t},\forall t\} \in \mathcal{F}_c,
\end{align}
where
\begin{align}
  ~ & F_{t}=(\pi_t^e + \pi_t^c \rho_t)(\hat{p}_{d,t}-\check{p}_{d,t}),\forall t,\label{eq:Ft}\\
  ~ & \mathcal{F}_c:=\{ \hat{p}_{d,t}, \check{p}_{d,t}, \forall t ~|~ \hat{p}_{d,t}=\sum_{v\in\mathcal{V}}\hat{p}^c_{v,t},\forall t, \label{equ:ubpdi}\\
  ~ & 0\leq \hat{p}_{v,t}^{c}\leq \bar{p}_{v},\forall v, \forall t\in[t_{v}^a,t_{v}^d],\label{equ:ubpd}\\
  ~ & \hat{e}_{v,t+1}= \hat{e}_{v,t}+\eta_c\hat{p}_{v,t}^{c}\Delta t,\forall v,\forall t \ne T,\label{equ:ube}\\
  ~ & \hat{e}_{v,t_v^a}= e_v^{a},~\hat{e}_{v,t_v^d}\geq e_v^{d},\forall v,\label{equ:ubsntatd}\\
  ~ & \underline{e}_v\leq \hat{e}_{v,t}\leq \bar{e}_v,\forall v,\forall t, \label{equ:uberanges}\\
  ~ & \check{p}_{d,t}=\sum_{v\in\mathcal{V}}\check{p}^c_{v,t},\forall t,\label{equ:lbpdi}\\
  ~ & 0\leq \check{p}_{v,t}^{c}\leq \bar{p}_{v},\forall v,\forall t\in[t_{v}^a,t_{v}^d],\\
  ~ & \check{e}_{v,t+1}= \check{e}_{v,t}+\eta_c\check{p}_{v,t}^{c}\Delta t,\forall v,\forall t \ne T,\label{equ:lbe}\\
  ~ & \check{e}_{v,t_v^a}= e_v^{a},~\check{e}_{v,t_v^d}\geq e_v^{d},\forall v,\\
  ~ & \underline{e}_v\leq \check{e}_{v,t}\leq \bar{e}_v,\forall v,\forall t,\label{equ:lberanges}\\
  ~ & \check p_{d,t} \leq \hat p_{d,t},\forall t, \label{equ:joint}\\
  ~ & \hat{p}^c_{v,t}=0, \check{p}^c_{v,t}=0, \forall v,\forall t \notin[t_{v}^a,t_{v}^d] \}. \label{eq:evava}
\end{align}
\end{subequations}
The goal of problem \textbf{P1} is to maximize the total value of aggregate EV power flexibility.
In \eqref{eq:p1-obj} and \eqref{eq:Ft}, $\pi_t,\forall t$ integrates the electricity price $\pi_t^e$ and carbon price $\pi_t^c$ to represent the unit value of EV power flexibility at different time. $\rho_t$ represents the carbon intensity of the grid power at time $t$.
The expectation in \eqref{eq:p1-obj} is taken with respect to uncertainties such as $\pi_t^e, \pi_t^c, \rho_t$ and EV charging tasks.
Constraint \eqref{equ:ubpdi} defines the upper aggregate charging power trajectory $\hat{p}_{d,t}$, which is the sum of the charging power bounds $\hat{p}_{v,t}^c$ of individual EVs.
Constraint \eqref{equ:ubpd} gives the charging power bounds of EV $v$, where $\bar{p}_{v}$ denotes the maximum power.
\eqref{equ:ube} describes the EV energy dynamics $\hat{e}_{v,t}$ whose bounds are defined in \eqref{equ:uberanges}.
EV's initial energy and charging requirements are specified in constraint \eqref{equ:ubsntatd}.
Similar restrictions apply to the lower trajectory of the aggregate EV power flexibility region, denoted by \eqref{equ:lbpdi}-\eqref{equ:lberanges}.
Inequality \eqref{equ:joint} ensures that $\{\hat p_{d,t},\forall t\}$ is above $\{\check p_{d,t},\forall t\}$.
Constraint \eqref{eq:evava} restricts the available charging time of EVs.

The model \eqref{eq:flexibility} can output an effective aggregate power flexibility region as stated in Proposition \ref{prop-1}, the proof of which can be found in Appendix \ref{appendix-A}.

\begin{proposition} \label{prop-1}
For aggregate EV charging power $\{p_{d,t},\forall t\}$ that satisfies $p_{d,t} \in [\check p_{d,t}, \hat p_{d,t}],\forall t$, there always exists a disaggregate, feasible charging strategy for individual EVs.    
\end{proposition}

\subsection{Stage 2: CS Energy Management and Carbon Trading}
\begin{figure}[htbp]
  \centering
  \includegraphics[width=0.4\textwidth]{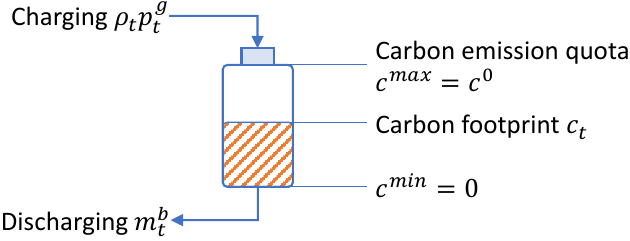}\\
  \caption{Generalized battery model for carbon footprint of EVCS.}\label{fig:carbonBat}
\end{figure}

Given the aggregate flexibility region by Stage 1, we can formulate the following offline problem to minimize the operation and carbon trading costs of the CS:
\begin{subequations}
\label{eq:evcsOperation}
\begin{align}
\textbf{P2:}~\min \limits_{p_{d,t},p^g_{t},m^b_{t},\forall t}  ~&\mathcal{C}:=\lim\limits_{T\rightarrow\infty}\frac{1}{T}\sum\limits_{t=1}^{T}\mathbb{E}\Big[C_{t}\Big],\label{eq:p2-obj}\\
  \text{s.t.}~&C_{t}=\pi^e_{t} p^g_{t} \Delta t + \pi^c_t m^b_{t},\forall t,\label{eq:Ct}\\
  &\check{p}_{d,t} \le p_{d,t} \le \hat{p}_{d,t}, \forall t, \label{eq:pidt-lbub} \\
  &p^g_{t} + p^r_{t} = p_{d,t}, \forall t, \label{eq:pgi}\\
  &0 \le p_t^r \le \bar{p}_t^r, \forall t, \label{eq:pri} \\
  &c_{t+1} = c_{t} + \rho_t p^g_{t} \Delta t - m^b_{t}, \forall t, \label{eq:c-state}\\
  &0 \le c_{t} \le c^0, \forall t, \label{eq:c-lbub}\\
  &0 \le m^b_{t} \le m^{b,max}, \forall t \in \mathcal{T_C}, \label{eq:cb-lbub}\\
  &m^b_{t} = 0, \forall t \notin \mathcal{T_C}. \label{eq:cb-0}
\end{align}
\end{subequations}
In \eqref{eq:p2-obj}, CSO aims to minimize the time average total cost, where $C_t$ includes two terms: the operation cost and carbon trading cost of CS as described in \eqref{eq:Ct}.
The expectation in \eqref{eq:p2-obj} is taken w.r.t uncertainties: electricity price $\pi_t^e$, carbon price $\pi_t^c$, grid carbon emission intensity $\rho_t$, and renewable generation $p_t^r$.
$p^g_t$ is the grid power purchased by the CS, and $m_{t}^b$ denotes the carbon emission quota purchased from the carbon market by CS.
Constraint \eqref{eq:pidt-lbub} indicates that the EV charging power should be within the aggregate flexibility region provided by Stage 1.
\eqref{eq:pgi} represents the power balance condition, where $p_t^r$ is the renewable generation power.
Constraint \eqref{eq:pri} represents that renewable power $p_t^r$ cannot exceed the expected maximum power $\bar{p}^r_{t}$.
Constraints \eqref{eq:c-state}-\eqref{eq:cb-lbub} present the generalized battery model for carbon footprint, as shown in Fig. \ref{fig:carbonBat}.
Particularly, \eqref{eq:c-state} represents the carbon footprint dynamics. $c_{t}$ is the carbon footprint of CS at time $t$, $\rho_t p_t^g \Delta t$ represents the produced carbon emission associated with the grid power consumption at time $t$, and $m^b_{t}$ denotes the carbon trading quantity of the CS at time $t$. $\rho_t p_t^g \Delta t$ increases the carbon footprint (similar to charging the battery) and $m^b_{t}$ reduces the carbon footprint (similar to discharging the battery).
Constraint \eqref{eq:c-lbub} describes that the CS's carbon footprint should always be within its carbon emission quota over the period.
Constraint \eqref{eq:cb-lbub} restricts the carbon trading quantity and trading time. $\mathcal{T_C}$ represents the set of time slots where carbon trading can take place.
Constraint \eqref{eq:cb-0} represents that in other time slots, there is no carbon trading.

The offline problems \textbf{P1} and \textbf{P2} effectively model how the two stages interact. However, they still face some solution challenges. That is, solving the offline problems \textbf{P1} and \textbf{P2} needs complete knowledge of future realization of uncertainties, such as future EV charging tasks, electricity prices, renewable generation power, carbon emission intensity, and carbon prices. However, these information are hardly available in practice. To address this challenge, online algorithms for the above two offline models are presented in the next sections.

\section{Online Characterization of Aggregate EV Power Flexibility}\label{sec:online-evagg}
In this section, we transform the offline problem \textbf{P1} into an online counterpart based on Lyapunov optimization.

\subsection{Problem Modification and EV Charging Queue Design}\label{subsec:pm}
Here, to accommodate the varying charging delays $t_v^d-t_v^a$ among different EVs, we collect and classify EVs' charging tasks into $G$ groups according to the charging delay. Each group is indexed by $g\in\mathcal{G}=\{1,2,...,G\}$ with delay $R_g$. $\hat{x}_{g,t}$ and $\check{x}_{g,t}$ are decision variables representing the upper and lower bound charging power for group $g$.
Therefore, corresponding to \textbf{P1}, we have $\hat{p}_{d,t}=\sum_g \hat{x}_{g,t}$ and $\check{p}_{d,t}=\sum_g \check{x}_{g,t}$.
Problem \textbf{P1} is then reformulated to evaluate the aggregate EV power flexibility characterization:
\begin{subequations}\label{eq:flexibility-p2}
\begin{align}
  \min\limits_{\hat{x}_{g,t},\check{x}_{g,t},\forall t} ~& \lim\limits_{T\rightarrow\infty}\frac{1}{T}\sum\limits_{t=1}^{T}\mathbb{E}\Big[-F_{t}\Big],\label{eq:P2''-obj} \\
  \text{s.t.}~& \lim\limits_{T\rightarrow\infty} \frac{1}{T}{\sum_{t=1}^T \mathbb{E}[\hat{a}_{g,t}-\hat x_{g,t}]} \le 0 ,\forall g, \label{eq:Qgub-lim} \\
  ~& \lim\limits_{T\rightarrow\infty}\frac{1}{T}{\sum_{t=1}^T \mathbb{E}[\check{a}_{g,t}-\check{x}_{g,t}]} \le 0,\forall g, \label{eq:Qglb-lim}\\
  ~& \hat{a}_{g,t} = \sum_{v\in \mathcal{V}_g}\hat{a}_{v,t},~
    \check{a}_{g,t} = \sum_{v \in \mathcal{V}_g}\check{a}_{v,t},\\
  ~& 0\leq\hat{x}_{g,t}\leq x_{g,max},\forall g,\forall t, \label{eq:xgub}\\
  ~& 0\leq\check{x}_{g,t}\leq x_{g,max},\forall g,\forall t, \label{eq:xglb}\\
  ~& \hat{x}_{g,t} \ge \check{x}_{g,t},\forall g,\forall t,  \label{eq:xgublb}
\end{align}
\end{subequations}
where $\hat{a}_{g,t}$ and $\check{a}_{g,t}$ are the upper and lower bounds of aggregate arrival charging demand of group $g$, determined by $\hat{a}_{g,t} = \sum_{v\in \mathcal{V}_g}\hat{a}_{v,t},~\check{a}_{g,t} = \sum_{v \in \mathcal{V}_g}\check{a}_{v,t}$.
The arrival desired charging demand of each EV $\check{a}_{v,t}$ and $\check{a}_{v,t}$ can be evaluated using the following equation, where each EV charging request is processed immediately upon arrival at the charging station:
\begin{align}
    & \check{a}_{v,t}=\left\{
    \begin{array}{ll}
    \bar{p}_{v},     & t_{v}^a\leq t<\lfloor \check{t}_{v}^{-}\rfloor+t_{v}^a \\
    \check{e}_{v}^{cha}/\eta_{c}-\lfloor \check{t}_{v}^{-}\rfloor \bar{p}_{v},  & t=\lfloor \check{t}_{v}^{-}\rfloor+t_{v}^a \\
    0,           & \text{otherwise}
    \end{array}\right. \nonumber\\ 
    & \hat{a}_{v,t}=\left\{
    \begin{array}{ll}
    \bar{p}_{v},     & t_{v}^a\leq t<\lfloor \hat{t}_{v}^{-}\rfloor+t_{v}^a \\
    \hat{e}_{v}^{cha}/\eta_{c}-\lfloor \hat{t}_{v}^{-}\rfloor \bar{p}_{v},     & t=\lfloor \hat{t}_{v}^{-}\rfloor+t_{v}^a \\
    0,           & \text{otherwise.}
    \end{array}\right. \nonumber 
\end{align}
Here, $\check{e}_{v}^{cha}=e_{v}^{d}-e_{v}^{a}$, and $\check{t}_{v}^{-}$ denotes the minimum needed charging time for EV $v$. It can be calculated by $\check{t}_{v}^{-}=\frac{\check{e}_{v}^{cha}}{\bar{p}_{v}\eta_{c}}$. $\lfloor.\rfloor$ means rounding down to the nearest integer.
Unlike $\check{a}_{v,t}$, we use the maximum charging demand $\bar{e}_{v}$ to calculate the upper bound of arrival charging demand $\hat{a}_{v,t}$.
In this case, we have $\hat{e}_{v}^{cha}=\bar{e}_{v}-e_{v}^{a}$, and $\hat{t}_{v}^{-}=\hat{e}_{v}^{cha}/(\bar{p}_{v}\eta_{c})$.

Constraint \eqref{eq:Qgub-lim} is a time-average constraint, meaning that the charging needs of all EVs can be satisfied by following the upper bound trajectory $\{\hat{x}_{g,t},\forall t\}$.
A similar condition applies to the lower bound trajectory $\{\check{x}_{g,t},\forall t\}$, as specified by constraint \eqref{eq:Qglb-lim}. 
The upper and lower trajectories of group $g$, $\{\hat{x}_{g,t},\check{x}_{g,t},\forall t\}$, are constrained by \eqref{eq:xgub} and \eqref{eq:xglb}.
Constraint \eqref{eq:xgublb} ensures the upper bound be no less than the lower bound.
$x_{g,max}=\sum_{v \in \mathcal{V}_g} \bar{p}_v$.

Next, to deal with the time average constraints \eqref{eq:Qgub-lim} and \eqref{eq:Qglb-lim}, we relax \eqref{eq:Qgub-lim} and \eqref{eq:Qglb-lim} by introducing two charging queues corresponding to the upper bound and lower bound charging power trajectories, respectively.
\begin{subequations}
\begin{gather}
    \hat{Q}_{g,t+1}=\max[\hat{Q}_{g,t}-\hat{x}_{g,t},0]+\hat{a}_{g,t},\label{eq:Qgt-ub}\\
    \check{Q}_{g,t+1}=\max[\check{Q}_{g,t}-\check{x}_{g,t},0]+\check{a}_{g,t}.\label{eq:Qgt-lb}
\end{gather}
\end{subequations}
According to \eqref{eq:Qgt-ub}, we have 
\begin{align}\label{eq-1}
    \hat Q_{g,t+1} -\hat a_{g,t} \ge \hat Q_{g,t} -\hat x_{g,t},\forall t.
\end{align}
Furthermore, we sum \eqref{eq-1} up across all $t$ and divide both sides by $T$ to obtain
\begin{align}
     \frac{\mathbb{E}[\hat{Q}_{g,T+1}]-\mathbb{E}[\hat{Q}_{g,1}]}{T} \ge \frac{\sum_{t=1}^T \mathbb{E}[\hat{a}_{g,t}-\hat x_{g,t}]}{T}.
\end{align}
As $\hat{Q}_{g,1}$ is 0, if queue $\hat{Q}_{g,t}$ is mean-rate-stable, i.e., $\lim \limits_{T\rightarrow\infty}{\mathbb{E}[\hat{Q}_{g,T+1}]}/{T}=0$, then constraint \eqref{eq:Qgub-lim} is satisfied. A similar condition holds for  $\check{Q}_{g,t}$. Therefore, we can convert the two constraints \eqref{eq:Qgub-lim} and \eqref{eq:Qglb-lim} into the mean-rate-stable conditions for queue $\hat Q_{g,t}$ and queue $\check Q_{g,t}$, respectively.

\subsection{Lyapunov Optimization}
With the constructed EV charging queues, we propose the online algorithm based on Lyapunov optimization.

\subsubsection{Lyapunov Function}
A concatenated vector of queues is defined: $\boldsymbol{\Theta}_t=(\boldsymbol{\hat{Q}}_{t},\boldsymbol{\check{Q}}_{t})$, where
\begin{subequations}
\begin{align}
\boldsymbol{\hat{Q}}_t=(\hat{Q}_{1,t},...,\hat{Q}_{G,t}),~\boldsymbol{\check{Q}}_t=(\check{Q}_{1,t},...,\check{Q}_{G,t}).
\end{align}
\end{subequations}
Then we define the Lyapunov function
\begin{equation}\label{equ:LyaFun}
L(\boldsymbol{\Theta}_t)=\frac{1}{2}\sum\limits_{g\in{\mathcal{G}}}\hat{Q}_{g,t}^2 + \frac{1}{2}\sum\limits_{g\in{\mathcal{G}}}\check{Q}_{g,t}^2,
\end{equation}
where $L(\boldsymbol{\Theta}_t)$ can be interpreted as a measure of the queue size. A smaller $L(\boldsymbol{\Theta}_t)$ is desirable because it indicates less congested queues $\hat{Q}_{g,t}$ and $\check{Q}_{g,t}$.

\subsubsection{Lyapunov Drift}
The conditional one-time slot Lyapunov drift is defined as follows:
\begin{equation}\label{equ:LyaDrift}
\Delta(\boldsymbol{\Theta}_t)=\mathbb{E}[L(\boldsymbol{\Theta}_{t+1}) - L(\boldsymbol{\Theta}_t)|\boldsymbol{\Theta}_t],
\end{equation}
where the expectation is conditional on the random factor $\boldsymbol{\Theta}_t$.

Given the present state $\boldsymbol{\Theta}_t$, $\Delta(\boldsymbol{\Theta}_t)$ measures the expected rise of the queue size.
It makes sense that reducing the Lyapunov drift would aid in virtual queue stabilization. However, concentrating only on reducing the Lyapunov drift could lead to a lower overall EV power flexibility value. 
To address this  issue, we add the objective function \eqref{eq:Ft} for time slot $t$ to \eqref{equ:LyaDrift}.
The drift-plus-penalty term is obtained,
\begin{equation}\label{equ:driftPlusP}
\Delta(\boldsymbol{\Theta}_t)+V_1\mathbb{E}[-F_t|\boldsymbol{\Theta}_t],
\end{equation}
where $V_1$ is a weight parameter that regulates the trade-off between maximizing aggregate EV power flexibility and queue stability.

\subsubsection{Minimizing Upper Bound} As \eqref{equ:driftPlusP} is still time-coupled due to $\Delta(\boldsymbol{\Theta}_t)$, we make further relaxation to replace \eqref{equ:driftPlusP} with its upper bound. Specifically, 
\begin{align}
&L(\boldsymbol{\Theta}_{t+1}) -L(\boldsymbol{\Theta}_t)\nonumber\\ &=\frac{1}{2}\sum\limits_{g\in{\mathcal{G}}}\Big\{\left[\hat{Q}_{g,t+1}^2-\hat{Q}_{g,t}^2\right] + \left[\check{Q}_{g,t+1}^2-\check{Q}_{g,t}^2\right]\Big\}.\label{equ:lyaD}
\end{align}
Using the queue $\hat{Q}_{g,t}$ update equation \eqref{eq:Qgt-ub}, we obtain
\begin{align}
\hat{Q}_{g,t+1}^2 &= \{\max[\hat Q_{g,t}-\hat x_{g,t},0]+\hat a_{g,t}\}^2\\
                  &\leq \hat Q_{g,t}^2+\hat{a}_{g,max}^2+\hat{x}_{g,max}^2+2\hat Q_{g,t}(\hat a_{g,t}-\hat x_{g,t}).\nonumber
\end{align}
Therefore,
\begin{align}
\frac{1}{2}[\hat{Q}_{g,t+1}^2-\hat{Q}_{g,t}^2]\leq \frac{1}{2}\left(\hat{x}_{g,max}^2+\hat{a}_{g,max}^2\right) \nonumber\\
+\hat{Q}_{g,t}\left(\hat{a}_{g,t}-\hat{x}_{g,t}\right).\label{eq:Q2-Ub}
\end{align}
Similarly, for queues $\check{Q}_{g,t}$,
\begin{align}
\frac{1}{2}\left[\check{Q}_{g,t+1}^2-\check{Q}_{g,t}^2\right]\leq \frac{1}{2}\left(\check{x}_{g,max}^2+\check{a}_{g,max}^2\right)\nonumber\\
+\check{Q}_{g,t}\left(\check{a}_{g,t}-\check{x}_{g,t}\right).\label{eq:Q2-lb}
\end{align}

Next, we replace the drift-plus-penalty term with the inequalities \eqref{eq:Q2-Ub} and \eqref{eq:Q2-lb}, yielding
\begin{equation}\label{equ:dppInequ}
\begin{aligned}
    &\Delta(\boldsymbol{\Theta}_t)+V_1\mathbb{E}[-F_t|\boldsymbol{\Theta}_t]\\
    &\leq A_1 + V_1\mathbb{E}[-F_t|\boldsymbol{\Theta}_t]+\sum\limits_{g\in{\mathcal{G}}}\hat{Q}_{g,t}\mathbb{E}\left[\hat{a}_{g,t}-\hat{x}_{g,t}|\boldsymbol{\Theta}_t\right]\\
    &+\sum\limits_{g\in{\mathcal{G}}}\check{Q}_{g,t}\mathbb{E}\left[\check{a}_{g,t}-\check{x}_{g,t}|\boldsymbol{\Theta}_t\right],
\end{aligned}
\end{equation}
where $A_1$ is a constant,
\begin{align*}
A_1 = ~ & \frac{1}{2}\sum\limits_{g\in{\mathcal{G}}}(\hat{x}_{g,max}^2+\hat{a}_{g,max}^2)+\frac{1}{2}\sum\limits_{g\in{\mathcal{G}}}(\check{x}_{g,max}^2+\check{a}_{g,max}^2).
\end{align*}

Finally, we can construct the following online optimization problem by rearranging the expression in \eqref{equ:dppInequ} and disregarding the constant terms:
\begin{align}
  \textbf{P1}': \min\limits_{\hat{x}_{g,t},\check{x}_{g,t}}~ &\sum\limits_{g\in{\mathcal{G}}}(-V_1\pi_t-\hat{Q}_{g,t})\hat{x}_{g,t}
+(V_1\pi_t-\check{Q}_{g,t})\check{x}_{g,t},\label{equ:P3}\\
  \hbox{s.t.}~ &\eqref{eq:xgub},\eqref{eq:xglb},\eqref{eq:xgublb},\nonumber
\end{align}
where, prior to solving $\textbf{P1}'$ in each time slot, $\hat{Q}_{g,t}$ and $\check{Q}_{g,t}$ are updated based on \eqref{eq:Qgt-ub} and \eqref{eq:Qgt-lb}.

The proposed approach solves problem $\textbf{P1}'$ to find the upper and lower bounds $\hat{x}_{g,t}$ and $\check{x}_{g,t}$ of the aggregate EV power flexibility region in the current time slot $t$, given the current system queue state $\boldsymbol{\Theta}_t$. As such, $\textbf{P1}$, the original offline optimization problem, has been separated into online problems.
The following statement can be used to calculate the difference between the online problem $\textbf{P1}'$ and the offline problem \textbf{P1}; the proof of this can be found in Appendix \ref{appendix-B}.
\begin{proposition}\label{prop-2}
Denote the obtained long-term time-average aggregate EV power flexibility value, specified in \eqref{eq:p1-obj}, of $\textbf{P1}$ and $\textbf{P1}'$ by $\mathcal{F}^{off}$ and $\mathcal{F}^{*}$, respectively. We have
\begin{equation}\label{eq:gap-1}
0 \le -\mathcal{F}^{*} + \mathcal{F}^{off} \leq \frac{1}{V_1}A_1,
\end{equation}
where $A_1$ is a constant defined in \eqref{equ:dppInequ}. 
\end{proposition}


\section{Online Algorithm for CS Energy Management and Carbon Trading}\label{sec:online-evcs}
In this section, constraint modification and virtual queue design are implemented to deal with the time-coupling constraints of carbon evolution to derive an online algorithm for Stage 2.

\subsection{Problem Modification}
We first convert the time-coupling constraint \eqref{eq:c-state} into a time-average one.
Both sides of \eqref{eq:c-state} are summed over $t\in\{1,...,T\}$ and then divided by $T$:
\begin{align}
\frac{1}{T}\sum_{t=1}^{T} m_{c,t}= \frac{c_{T+1}}{T}-\frac{c_{1}}{T}, \label{eq:cAvg1}
\end{align}
where $m_{c,t}=-m^b_{t}+\rho_t p^{g}_{t}$.
We take expectations on both sides of \eqref{eq:cAvg1} and then let $T$ go to infinity:
\begin{equation}\label{eq:cAvg2}
\lim\limits_{T\rightarrow\infty}\frac{1}{T}\sum\limits_{t=1}^{T}\mathbb{E}\left[m_{c,t}\right]= \lim\limits_{T\rightarrow\infty}\mathbb{E}\left[\frac{c_{T+1}}{T}\right]-\lim\limits_{T\rightarrow\infty}\mathbb{E}\left[\frac{c_{1}}{T}\right],
\end{equation}

Due to the \eqref{eq:c-lbub}, $c_{1}$ and $c_{T+1}$ are finite, hence the right side of \eqref{eq:cAvg2} equals zero. As a result,
\begin{equation}\label{eq:cAvg3}
\lim\limits_{T\rightarrow\infty}\frac{1}{T}\sum\limits_{t=1}^{T}\mathbb{E}\left[m_{c,t}\right]= 0.
\end{equation}
Note that constraint \eqref{eq:cAvg3} is a relaxed version of constraints \eqref{eq:c-state}-\eqref{eq:c-lbub}. The above relaxation step facilitates the implementation of Lyapunov optimization techniques.


\subsection{Virtual Carbon Queue Design}
The virtual carbon queue $H_{t}$ is defined as follows:
\begin{equation}\label{eq:H}
H_{t}= c_{t}-\phi_{t},
\end{equation}
where $\phi_{t}$ is a perturbation parameter designed to ensure the feasibility of constraint \eqref{eq:c-lbub}, which is explained later.
Then we can write the dynamics of virtual carbon queue:
\begin{equation}\label{eq:HDyn}
H_{t+1}= H_{t}+m_{c,t}.
\end{equation}
By comparing \eqref{eq:HDyn} with \eqref{eq:c-state}, we can observe that $H_{t}$ is a shifted version of $c_t$. But different from $c_{t}$, the virtual carbon queue $H_{t}$ can take on negative values due to the perturbation parameter $\phi_{t}$.
This shift ensures that the constraint \eqref{eq:c-lbub} is satisfied. Furthermore, due to \eqref{eq:cAvg3} and $c_1$ is zero, the virtual carbon queue $H_t$ satisfies the mean-rate-stable condition, i.e.,
$\lim\limits_{t\rightarrow\infty}{\mathbb{E}[H_t]}/{t}=0.$

\subsection{Lyapunov Optimization} \label{sec:III-B}

\subsubsection{Lyapunov Function}
The Lyapunov function of queue $H_{t}$ is given by
\begin{equation}\label{eq:LyaFun}
L(H_{t})=\frac{1}{2}H_{t}^2,
\end{equation}

\subsubsection{Lyapunov Drift of $L(H_{t})$}
\begin{equation}\label{eq:LyaDrift}
\Delta(H_{t})=\mathbb{E}[L(H_{t+1}) - L(H_{t})|H_{t}],
\end{equation}
where the expectation is conditional on the random $H_{t}$.

\subsubsection{Drift-Plus-Penalty}
The drift-plus-penalty term is obtained as follows
\begin{equation}\label{eq:driftPlusP}
\Delta(H_{t})+V_2\mathbb{E}[C_{t}|H_{t}],
\end{equation}
where the weight parameter $V_2$ manages the trade-off between the stability of virtual queues and the reduction of operation cost. We provide Proposition \ref{prop-2} later to explain how to determine the value of $V_2$.

\subsubsection{Minimizing Upper Bound}
Due to $\Delta(H_{t})$, problem \eqref{eq:driftPlusP} remains time-coupled. We minimize its upper bound to derive the control decision, which can be adapted to an online implementation. First, the Lyapunov drift for one time slot is obtained:
\begin{align}
L(H_{t+1}) -L(H_{t}) =&\frac{1}{2}\left[H_{i,t+1}^2-H_{t}^2\right].\label{eq:lyaD}
\end{align}
Based on the queue update equation \eqref{eq:HDyn}, we have
\begin{align}
\frac{1}{2}&\left[H_{i,t+1}^2-H_{t}^2\right] 
\le H_{t}m_{c,t} + \frac{1}{2}\max\{(\rho_t p^g_{t})^2,(m^{b,max})^2\},\label{eq:H-ub}
\end{align}
Then the drift-plus-penalty term satisfies
\begin{equation}\label{eq:dppInequ}
\begin{aligned}
&\Delta(H_t)+V_2\mathbb{E}[C_t|H_t] 
\leq A_2 + H_{t}\mathbb{E}\left[m_{c,t}|H_t\right] + V_2\mathbb{E}[C_t|H_t],
\end{aligned}
\end{equation}
where $A_2 = \frac{1}{2}\max\{(\rho^{max}_t p^{g,max}_{t})^2,(m_{t}^{b,max})^2\}$.

By minimizing the upper bound of the drift-plus-penalty term, we obtain the following online optimization problem
\begin{subequations}
\begin{align}
  \textbf{P2}': \min ~ &H_{t}m_{c,t} + V_2 C_t,\label{eq:P3-2}\\
  \hbox{s.t.}~ &\eqref{eq:pidt-lbub},\eqref{eq:pgi},\eqref{eq:pri},\eqref{eq:cb-lbub}, \\
  &\check{p}_{d,t}=\sum_{g\in\mathcal{G}} \check{x}_{g,t}, ~ \hat{p}_{d,t}=\sum_{g\in\mathcal{G}} \hat{x}_{g,t},
\end{align}
\end{subequations}
By solving $\textbf{P2}'$, the proposed method determines $p_{d,t},p_t^g,m_t^b$ in each time slot $t$ given the virtual queue state $H_t$.

\subsection{Analysis of Feasibility and Performance}
Comparing constraints of $\textbf{P2}$ with those of $\textbf{P2}'$, constraint \eqref{eq:c-lbub} is not explicitly considered in $\textbf{P2}'$.
In fact, the bound constraint \eqref{eq:c-lbub} of carbon footprint can be guaranteed by carefully choosing the perturbation parameter $\phi_{t}$, as indicated in the following proposition.
\begin{proposition} \label{prop-3}
When $\rho_t^{max} p^{g,max}_{t} + m^{b,max} \le c^{max}$ holds, if we let
\begin{equation}\label{eq:phi}
\phi_{t} = m^{b,max}+V_2\pi_t^{g,max}\frac{1}{\rho_t}, \forall t,
\end{equation}
where
\begin{align}
0\leq V_2 \leq V_{2,max} = \frac{c^{max}-m^{b,max}}
                      {\pi_{t}^{c,max}+\pi_{t}^{e,max}\frac{1}{\rho_t^{min}}}\label{eq:Vmax},
\end{align}
then the sequence of optimal solutions obtained by online problem $\textbf{P2}'$ can satisfy the constraint \eqref{eq:c-lbub}.
\end{proposition}

Proposition \ref{prop-3} is proven in Appendix \ref{appendix-C}. Furthermore, the difference in optimal solutions between the offline problem \textbf{P2} and the online problem $\textbf{P2}'$ is discussed below.
\begin{proposition}
\label{prop-4}
If let $C^*$ and $\widehat{C}$ represent the attained long-term time-average cost objective values of $\textbf{P2}$ and $\textbf{P2}'$, respectively, we have
\begin{equation}\label{eq:gap-2}
0 \le \widehat{C} - C^* \leq \frac{1}{V_2}A_2,
\end{equation}
with constant $A_2$ defined in \eqref{eq:dppInequ}.
\end{proposition}

Please refer to Appendix \ref{appendix-D} for the proof of Proposition \ref{prop-4}.
Parameter $V_2$ influences the optimality gap.
A larger $V_2$ value leads to a large virtual queue, but it can narrow the optimality gap.
Conversely, a smaller $V_2$ value leads to more stable queues, but it also widens the optimality gap.

\subsection{Two-stage Online Algorithm}
Algorithm \ref{algo} provides a detailed description of the proposed two-stage online algorithm for CS operation, which includes aggregate EV power flexibility characterization, energy management, and carbon trading.

\begin{algorithm}[htbp]
\caption{Two-stage online algorithm}\label{algo}
\begin{algorithmic}[1]
\renewcommand{\algorithmicrequire}{ \textit{\textbf{Initialization}}}
\REQUIRE
\STATE{Set $t=1$, queues $\hat{Q}_{g,1}=0$, $\check{Q}_{g,1}=0$, $H_1=0$, parameters $V_1>0,V_2>0$.}
\renewcommand{\algorithmicrequire}{ \textit{\textbf{Stage 1: EV Aggregation}}}
\REQUIRE
\STATE{CSO classifies the arriving EVs and places them into different queues $\hat{Q}_{g,t}$, $\check{Q}_{g,t}$ based on $R_g$.}
\STATE{\color{black}Aggregate EV power flexibility region $[\check{p}_{d,t}^*,\hat{p}_{d,t}^*]$ is obtained by solving problem $\textbf{P1}'$.}
\renewcommand{\algorithmicrequire}{ \textit{\textbf{Stage 2: Energy Management and Carbon Trading}}}
\REQUIRE
\STATE{Solve problem $\textbf{P2}'$ to determine the dispatch power of EVs, grid power, renewable power, and carbon trading quantity.}
\STATE{Update queue $c_{t+1}$ and queue $H_{t+1}$.}
\STATE{Update queues $\hat{Q}_{g,t+1}$, $\check{Q}_{g,t+1}$ based on \eqref{eq:Qgt-ub}, \eqref{eq:Qgt-lb}.}
\STATE{$t \leftarrow t+1$, and repeat Stages 1-2.}
\end{algorithmic}
\end{algorithm}

\section{Case Studies}
\label{secV}
In this section, we evaluate the performance of the proposed two-stage online algorithm on a CS located in a commercial area and compare it with other traditional methods.

\subsection{System Setup} \label{sec:case-A}
\begin{figure}[!htbp]
  \centering
  \includegraphics[width=0.4\textwidth]{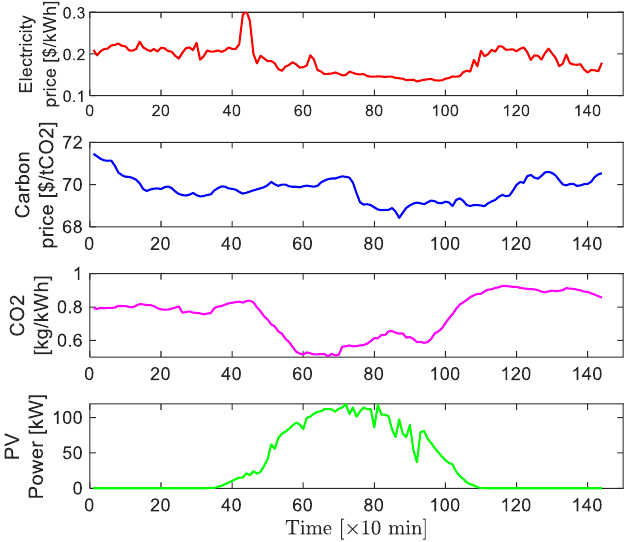}\\
  \caption{Real-time electricity price, carbon trading price, carbon emission intensity, and PV power generation profiles.}\label{fig:price}
\end{figure}

\begin{figure}[!htbp]
  \centering
  \includegraphics[width=0.4\textwidth]{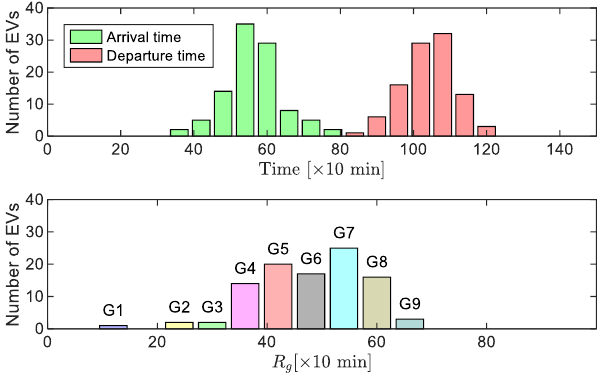}\\
  \caption{Distribution of EV arrival/departure times, and charging delay $R_g$.}\label{fig:time}
\end{figure}

The time period considered in the simulations is 24 hours, divided into 144 time slots, i.e. each time slot covers 10 minutes.
Carbon trading time is set as $\mathcal{T_C}=[6,12,18,\ldots,144]$ rather than all time slots.
Fig. \ref{fig:price} shows the dynamic real-time data of electricity price, carbon trading price, carbon emission intensity, and PV power generation \cite{pjm}, respectively.
As seen from the figure, those data profiles exhibit high uncertainty.
A total of 100 EVs located within a commercial area are considered. Suppose their arrival and departure times follow Gaussian distributions: $t_v^a \in \mathcal{N}(\text{9:00}, (1.2hr)^2)$ and $t_v^a \in \mathcal{N}(\text{18:00}, (1.2hr)^2)$ \cite{mohamed2014real}, as illustrated in Fig. \ref{fig:time} (top). The parking duration $R_g$ ranges from 2 to 11 hours, as illustrated in Fig. \ref{fig:time} (bottom). There are $G=9$ groups of EVs. The initial battery energy level of each EV is uniformly distributed within $[0.3,0.5] \times$ the EV battery capacity \cite{jin2014optimized}. The required SoC upon departure is set to be 0.5, and the maximum SoC upon departure is no greater than 0.9.
In addition, three typical EV battery capacities-24, 40, and 60 kWh-as well as three charging power levels-3.3, 6.6, and 10 kW-are taken into consideration to reflect the diversity of EVs.
The weight parameter $V_1$ is set as 20, and $V_2$ satisfies \eqref{eq:Vmax}.
The initial carbon emission quota of CS is 80 kgCO2.

\subsection{Results Analysis}

\begin{figure}[!htbp]
  \centering
  \includegraphics[width=0.4\textwidth]{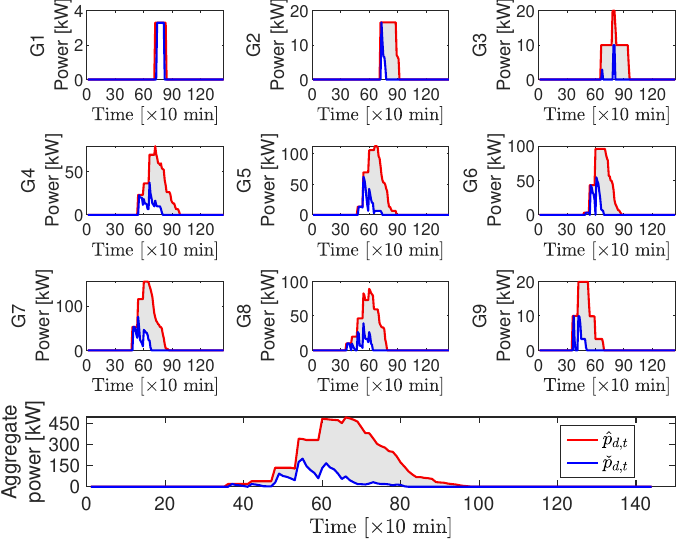}\\
  \caption{The obtained aggregate EV power flexibility region of each group.}\label{fig:agg}
\end{figure}

\begin{figure}[!htbp]
  \centering
  \includegraphics[width=0.4\textwidth]{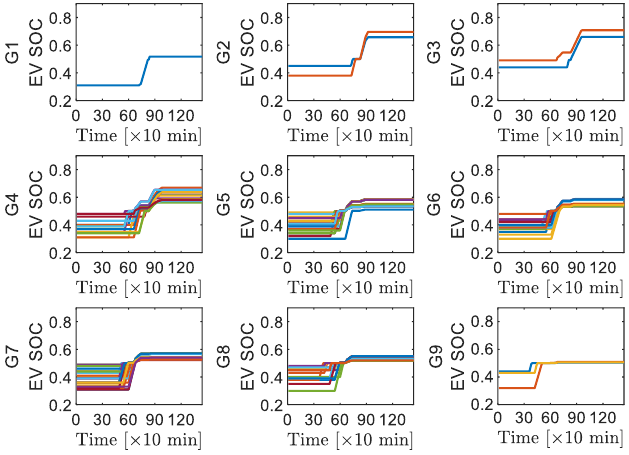}\\
  \caption{EV SoC profile of each group.}\label{fig:evsoc}
\end{figure}

We first investigate the EV aggregation. By running Algorithm \ref{algo}, we can obtain the optimal upper and lower aggregate power boundaries, $\hat p_{d,t}$ and $\check p_{d,t}$, the area between which is the EV aggregate power flexibility region.
Fig. \ref{fig:agg} illustrates the aggregate EV power flexibility regions for each group, as well as the total region.
First, we can find that the optimized lower aggregate power boundary $\check p_{d,t}$ is not constantly zero. While setting the EV charging power to zero at a single time slot may maximize the power flexibility region, an all-zero strategy will fail to meet the energy requirement when EVs depart. 
This is due to the time-coupling nature of EV charging.
The proposed algorithm effectively tackles this issue, as evidenced by the determined non-all-zero $\check p_{d,t}$, which ensures that EVs can reach the desired SoC upon departure. As depicted in Fig. \ref{fig:evsoc}, the disaggregation and dispatch results demonstrate that all EVs attain a higher SoC than the targeted 0.5.
Furthermore, we can find that the total power flexibility region initially exhibits a relatively narrow range, gradually expanding over time. This arises from the fact that most EVs start charging earlier until they almost fulfill their charging requirements. Once they complete charging, the minimal charging power could be zero to maximize the available power flexibility. 
In addition, the boundaries are influenced by the parameter $V_1$. A larger $V_1$ emphasizes power flexibility, which will be discussed in detail later.

\subsection{Analysis of Energy Management and Carbon Trading}
\begin{figure}[!htbp]
  \centering
  \includegraphics[width=0.4\textwidth]{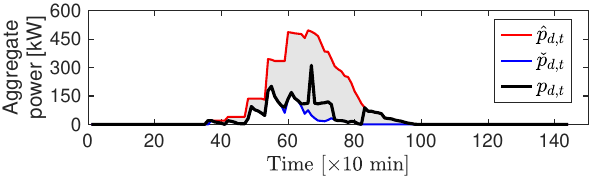}\\
  \caption{The aggregate dispatch profile in the second stage.}\label{fig:dispatchAgg}
\end{figure}

With the determined boundaries of aggregate EV power flexibility at each time slot, the CSO only needs to dispatch aggregate EV power instead of controlling EVs individually. This facilitates the energy management and carbon trading problem.
As shown in Fig. \ref{fig:dispatchAgg}, the optimal aggregate EV power dispatch remains close to the lower bound for most of the time, with a large margin from the upper bound $\hat p_{d,t}$.
On the one hand, this indicates that the available ramping-up aggregate EV power flexibility is sufficient. On the other hand, it implies that this strategy of lower power dispatch benefits in reducing grid power consumption, thereby lowering electricity costs and mitigating carbon emissions. Specifically, during the peak price period $t\in(43,45)$, the energy management strategy strategically minimizes the aggregate power dispatch and shifts the charging demand to later time slots with lower electricity prices. In contrast, larger EV charging powers are dispatched for $t\in (64,73)$ and $(81,97)$ when electricity prices are lower, to improve the economic efficiency.

\begin{figure}[!htbp]
  \centering
  \includegraphics[width=0.4\textwidth]{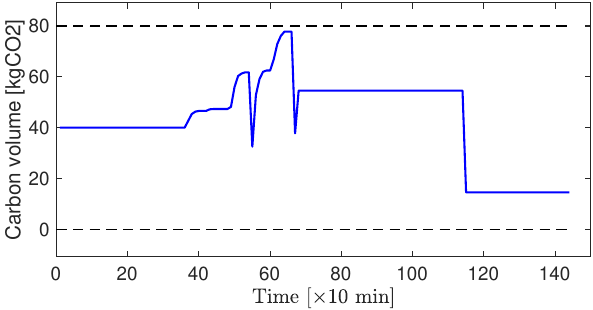}\\
  \caption{carbon footprint evolution over time.}\label{fig:cVol}
\end{figure}

Proposition \ref{prop-3} claims that by properly designing the parameter $\phi_t$, the carbon footprint bound constraint can be satisfied without explicitly being considered in $\textbf{P2}'$. Here, we present the related simulation results to verify the feasibility of the proposed algorithm.
Fig. \ref{fig:cVol} shows the evolution of carbon footprint over time, with an initial value of 40 kgCO2. First, the $c_t$ gradually increases due to EV charging and grid power consumption.
In the $t=54$ period, when the carbon footprint reaches 62, close to the carbon emission quota limit, carbon trading is triggered, resulting in a significant drop in the next time slot.
Afterward, the carbon footprint continues to grow until $t=66$, leading to the occurrence of the second carbon trade.
During $t\in(68,114)$, the carbon footprint remains unchanged, indicating no grid energy consumption. This is attributed to the renewable PV power generation supplies all the EV charging demand during this period.
At $t=114$, the third carbon trade takes place, causing a carbon footprint decline.
Following this, carbon footprint remains unchanged as there is no EV charging and grid power consumption.
Overall, the carbon footprint remains below the carbon emission quota throughout the entire process. This demonstrates the effectiveness of the proposed algorithm and Proposition \ref{prop-3}.
In addition, carbon trading occurs infrequently, aligning with the practice in the carbon trading market.

\subsection{Impact of Parameters}
\begin{figure}[!htbp]
  \centering
  \includegraphics[width=0.4\textwidth]{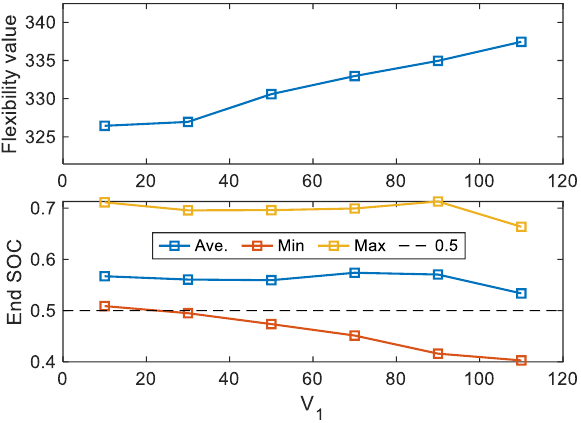}\\
  \caption{Impact of $V_1$ on the flexibility value and end SOC.}\label{fig:V1}
\end{figure}

As mentioned in \eqref{eq:driftPlusP}, parameter $V_1$ plays a crucial role in balancing the stability of EV charging queues and maximizing the total power flexibility value. Here, we change the value of $V_1$.
As shown in the top of Fig. \ref{fig:V1}, the total power flexibility value increases as $V_1$ increases. This is due to the increased emphasis placed on maximizing power flexibility.
The bottom figure summarizes the minimal, average, and maximal end SoC values of all EVs. Overall, the minimal end SoC of all EVs shows a decreasing trend as $V_1$ increases. Moreover, it gradually falls below the desired SoC value of 0.5 (dash line) when departure. This suggests that while a higher value of $V_1$ may provide a higher power flexibility value, it may also introduce instability to the charging queue, raising the risk of failing to meet the EV charging requirements.

\begin{figure}[!htbp]
  \centering
  \includegraphics[width=0.4\textwidth]{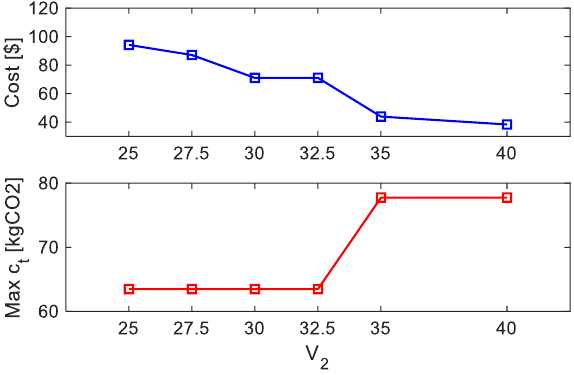}\\
  \caption{Impact of $V_2$ on the cost and maximal value of $c_t$.}\label{fig:V2}
\end{figure}

\begin{figure}[!htbp]
  \centering
  \includegraphics[width=0.4\textwidth]{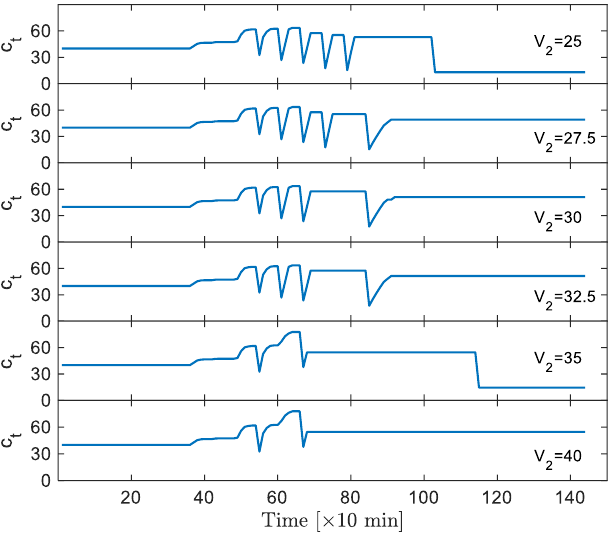}\\
  \caption{carbon footprint evolution $c_t$ over time under different $V_2$.}\label{fig:ct}
\end{figure}

Similarly, $V_2$ manages the trade-off between stabilizing the virtual carbon queue and minimizing the total cost of the CS.
Overall, the top of Fig. \ref{fig:V2} shows that the CS cost decreases with increasing $V_2$. This is because a larger $V_2$ places more emphasis on reducing the cost.
In contrast, the maximal carbon footprint $c_t$ recorded during the CS operation gradually increases as $V_2$ increases. A higher $c_t$ indicates it is close to the allowed carbon emission quota limit, i.e., the risk of instability for the virtual carbon queue. Conversely, a smaller $V_2$ places more emphasis on the stability of the carbon queue, resulting in a smaller carbon footprint.
Besides, Fig. \ref{fig:ct} presents the evolution of carbon footprint over time with different $V_2$. When $V_2$ is smaller, such as $V_2=10$, more carbon trading events occur, corresponding to the drop in carbon footprint. This is because preventing carbon footprint growth helps maintain the stability of virtual carbon queues. On the contrary, as $V_2$ increases, the number of carbon trading decreases, and the carbon footprint grows.
In addition, the carbon footprint still remains within the carbon emission quota under different values of $V_2$. This again verifies Proposition \ref{prop-3}.

\begin{figure}[!htbp]
  \centering
  \includegraphics[width=0.4\textwidth]{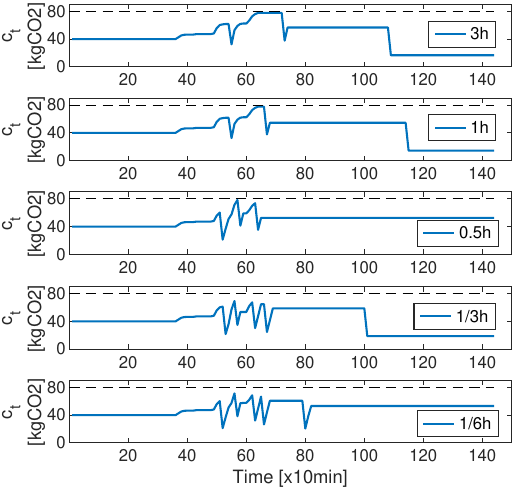}\\
  \caption{$c_t$ over time under different $\Delta t_c$.}\label{fig:ct-tc}
\end{figure}

\begin{figure}[!htbp]
  \centering
  \includegraphics[width=0.4\textwidth]{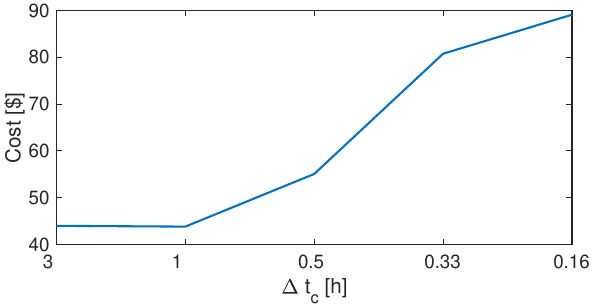}\\
  \caption{Total cost under different $\Delta t_c$.}\label{fig:cost-tc}
\end{figure}

We further analyze the impact of carbon trading time interval, $\Delta t_c$, on the results.
Fig. \ref{fig:ct-tc} illustrates the carbon footprint over time under different $\Delta t_c$. As $\Delta t_c$ decreases, more drops appear in the profile of $c_t$, indicating a higher frequency of carbon trading. This is because a smaller carbon trading time interval gives rise to more opportunities for carbon trading, allowing the system to maintain the stability of the carbon queue and not exceed the carbon emission quota. As a result, the total cost correspondingly increases, as shown in Fig. \ref{fig:cost-tc}.

\subsection{Performance Comparison}
The proposed online algorithm is compared with other benchmarks to demonstrate its advantage in terms of the total aggregate power flexibility value:
\begin{itemize}
\item Benchmark 1 (B1) is a charging requirement prioritized method that ensures EV charging requirements are satisfied first. EVs charge at maximum power upon arrival until reaching the desired SOC. Then, it focuses on maximizing power flexibility, allowing for flexible charging within SoC limits.
\item Benchmark 2 (B2) is an offline method that assumes perfect knowledge of future information. It directly solves problem \textbf{P1} to obtain the aggregate EV power flexibility region over the entire day. While it is unrealistic in practice, B2 serves as a theoretical benchmark to evaluate the performance of others.
\item Proposed online method with different dispatch ratios $\alpha$. 
\end{itemize}

Table \ref{tab:agg} presents the total flexibility value of different methods. B1 achieves the lowest total flexibility value, indicating that prioritizing charging requirements limits the overall power flexibility. The proposed method shows a significant increase in total flexibility value as the dispatch ratio increases. This is because a higher dispatch ratio allocates more charging energy to EVs, which alleviates EVs' charging anxiety and releases more flexibility. B2 achieves the highest power flexibility value but relies on perfect knowledge of future uncertainties, which is unrealistic in practice. Overall, the proposed method outperforms B1 and is more practical than B2.

\begin{table}[!htbp]
  \centering
  \caption{Total aggregate flexibility value comparison (Unit: USD).}\label{tab:agg}
  \begin{tabular}{ccccccc}
  \hline\hline
    Methods   & B1    & B2     & Prop.   & Prop.    & Prop.   & Prop. \\
    $\alpha$  & -     & -      & 0.2     & 0.4      & 0.6     & 0.8 \\
    \hline
    Value & 326.5     & 351.8  & 326.7   & 326.7    & 326.7   & 326.7\\
    Improvement & -   & 7.7\%  & 0.9\%   & 2.1\%    & 3.0\%   & 3.4\% \\
    \hline\hline
   \end{tabular}
\end{table}

Next, we compare the proposed method with other methods in terms of the total operation cost of the CS.
\begin{itemize}
    \item Benchmark 3 (B3) greedily minimizes the total cost in the current time slot and ignores the long-term benefits. In addition, carbon trading is executed when the carbon footprint reaches the carbon emission quota limit.
    \item Benchmark 4 (B4) is an offline method that assumes known future information of uncertainties.
\end{itemize}

Table \ref{tab:cost} lists the total cost obtained by different methods.
In the simulation of B3, operational infeasibility occurs when the carbon footprint approaches the carbon emission quota limit due to the inherent conflict between grid power consumption and the carbon emission quota constraint, which are coupled within the carbon dynamics \eqref{eq:c-state}. Specifically, a larger grid power consumption is required to meet the charging demand but it will violate the carbon emission quota limit, and vice versa. To make the problem feasible, we modified B3 by adding a penalty term for additional power consumption. In contrast, both the offline method and the proposed method are feasible. The offline method achieves the lowest cost, but it relies on perfect knowledge of future information. Overall, the proposed method achieves the trade-off between feasibility and practicality.

\begin{table}[!htbp]
  \centering
  \caption{Total cost comparison (Unit: USD).}\label{tab:cost}
  \begin{tabular}{ccccccc}
  \hline\hline
    Methods   & B3          & B3 modified   & B4    & Prop. \\
    \hline
    Value     & Infeasible  & 64.9          & 21.4  & 35.5 \\
    Reduction & -           & -             & 67\%  & 45\% \\
    \hline\hline
   \end{tabular}
\end{table}

\section{Conclusion}
\label{secVI}
This paper proposes a two-stage online algorithm to mitigate the increasing complexity of CS charging management by massive EVs and the substantial carbon emissions for CS. In the first stage, the aggregate EV charging power flexibility is characterized by lower and upper charging power bounds. It enables the CSO to treat all EVs as a whole for subsequent dispatch. In the second stage, we focus on energy management and carbon trading issues. A generalized battery model is proposed to capture the dynamics of the carbon emission level and carbon trading, adapting it to online carbon trading approaches. Charging queues and virtual carbon queues are designed, and Lyapunov optimization is employed in each stage to transform the offline models into their online counterparts. Case studies validate the effectiveness and benefits of the proposed method. We have the following findings:

1) The proposed online algorithm provides a near-optimal result in terms of flexibility and operation cost. The EV charging requirement is also satisfied before departure.

2) A larger parameter $V_1$ and $V_2$ lead to a higher total flexibility and a lower total cost but increase the risk of failing to meet the EV charging requirements and exceeding the carbon emission quota limit, respectively.

3) Carbon trading time interval influences the carbon trading number and total cost. A smaller carbon trading time interval leads to less carbon trading and lower costs.

\bibliographystyle{IEEEtran}
\bibliography{PaperRef}

\begin{thebibliography}{10}
\providecommand{\url}[1]{#1}
\csname url@samestyle\endcsname
\providecommand{\newblock}{\relax}
\providecommand{\bibinfo}[2]{#2}
\providecommand{\BIBentrySTDinterwordspacing}{\spaceskip=0pt\relax}
\providecommand{\BIBentryALTinterwordstretchfactor}{4}
\providecommand{\BIBentryALTinterwordspacing}{\spaceskip=\fontdimen2\font plus
\BIBentryALTinterwordstretchfactor\fontdimen3\font minus
  \fontdimen4\font\relax}
\providecommand{\BIBforeignlanguage}[2]{{%
\expandafter\ifx\csname l@#1\endcsname\relax
\typeout{** WARNING: IEEEtran.bst: No hyphenation pattern has been}%
\typeout{** loaded for the language `#1'. Using the pattern for}%
\typeout{** the default language instead.}%
\else
\language=\csname l@#1\endcsname
\fi
#2}}
\providecommand{\BIBdecl}{\relax}
\BIBdecl

\bibitem{dixon2020scheduling}
J.~Dixon, W.~Bukhsh, C.~Edmunds, and K.~Bell, ``Scheduling electric vehicle
  charging to minimise carbon emissions and wind curtailment,'' \emph{Renewable
  Energy}, vol. 161, pp. 1072--1091, 2020.

\bibitem{mouli2016system}
G.~C. Mouli, P.~Bauer, and M.~Zeman, ``System design for a solar powered
  electric vehicle charging station for workplaces,'' \emph{Appl. Energy}, vol.
  168, pp. 434--443, 2016.

\bibitem{yang2018distributed}
Y.~Yang, Q.-S. Jia, G.~Deconinck, X.~Guan, Z.~Qiu, and Z.~Hu, ``Distributed
  coordination of {EV} charging with renewable energy in a microgrid of
  buildings,'' \emph{IEEE Trans. Smart Grid}, vol.~9, no.~6, pp. 6253--6264,
  2018.

\bibitem{wu2017optimal}
X.~Wu, X.~Hu, Y.~Teng, S.~Qian, and R.~Cheng, ``Optimal integration of a hybrid
  solar-battery power source into smart home nanogrid with plug-in electric
  vehicle,'' \emph{J. Power Sources}, vol. 363, pp. 277--283, 2017.

\bibitem{wu2018stochastic}
X.~Wu, X.~Hu, X.~Yin, and S.~J. Moura, ``Stochastic optimal energy management
  of smart home with pev energy storage,'' \emph{IEEE Trans. Smart Grid},
  vol.~9, no.~3, pp. 2065--2075, 2018.

\bibitem{yan2018optimized}
Q.~Yan, B.~Zhang, and M.~Kezunovic, ``Optimized operational cost reduction for
  an {EV} charging station integrated with battery energy storage and {PV}
  generation,'' \emph{IEEE Trans. Smart Grid}, vol.~10, no.~2, pp. 2096--2106,
  2018.

\bibitem{xu2016hierarchical}
Z.~Xu, D.~S. Callaway, Z.~Hu, and Y.~Song, ``Hierarchical coordination of
  heterogeneous flexible loads,'' \emph{IEEE Trans. Power Syst.}, vol.~31,
  no.~6, pp. 4206--4216, 2016.

\bibitem{zhang2016evaluation}
H.~Zhang, Z.~Hu, Z.~Xu, and Y.~Song, ``Evaluation of achievable vehicle-to-grid
  capacity using aggregate {PEV} model,'' \emph{IEEE Trans. Power Syst.},
  vol.~32, no.~1, pp. 784--794, 2016.

\bibitem{wang2022evaluation}
L.~Wang, J.~Kwon, N.~Schulz, and Z.~Zhou, ``Evaluation of aggregated {EV}
  flexibility with tso-dso coordination,'' \emph{IEEE Trans. Sustain. Energy},
  vol.~13, no.~4, pp. 2304--2315, 2022.

\bibitem{shi2021optimal}
X.~Shi, Y.~Xu, Q.~Guo, and H.~Sun, ``Optimal dispatch based on aggregated
  operation region of {EV} considering spatio-temporal distribution,''
  \emph{IEEE Trans. Sustain. Energy}, vol.~13, no.~2, pp. 715--731, 2021.

\bibitem{zhou2021forming}
M.~Zhou, Z.~Wu, J.~Wang, and G.~Li, ``Forming dispatchable region of electric
  vehicle aggregation in microgrid bidding,'' \emph{IEEE Trans. Ind.
  Informat.}, vol.~17, no.~7, pp. 4755--4765, 2021.

\bibitem{chen2020aggregate}
X.~Chen, E.~Dall'Anese, C.~Zhao, and N.~Li, ``Aggregate power flexibility in
  unbalanced distribution systems,'' \emph{IEEE Trans. Smart Grid}, vol.~11,
  no.~1, pp. 258--269, 2020.

\bibitem{feng2024dayahead}
S.~Feng, W.~Wei, and Y.~Chen, ``Day-ahead scheduling and online dispatch of
  energy hubs: A flexibility envelope approach,'' \emph{IEEE Trans. Smart
  Grid}, vol.~15, no.~3, pp. 2723--2737, 2024.

\bibitem{wang2021aggregate}
S.~Wang and W.~Wu, ``Aggregate flexibility of virtual power plants with
  temporal coupling constraints,'' \emph{IEEE Trans. Smart Grid}, vol.~12,
  no.~6, pp. 5043--5051, 2021.

\bibitem{tulpule2013economic}
P.~J. Tulpule, V.~Marano, S.~Yurkovich, and G.~Rizzoni, ``Economic and
  environmental impacts of a pv powered workplace parking garage charging
  station,'' \emph{Appl. Energy}, vol. 108, pp. 323--332, 2013.

\bibitem{hoehne2016optimizing}
C.~G. Hoehne and M.~V. Chester, ``Optimizing plug-in electric vehicle and
  vehicle-to-grid charge scheduling to minimize carbon emissions,''
  \emph{Energy}, vol. 115, pp. 646--657, 2016.

\bibitem{li2023smart}
J.~Li, G.~Wang, X.~Wang, and Y.~Du, ``Smart charging strategy for electric
  vehicles based on marginal carbon emission factors and time-of-use price,''
  \emph{Sustainable Cities and Society}, vol.~96, p. 104708, 2023.

\bibitem{yan2021blockchain}
M.~Yan, M.~Shahidehpour, A.~Alabdulwahab, A.~Abusorrah, N.~Gurung, H.~Zheng,
  O.~Ogunnubi, A.~Vukojevic, and E.~A. Paaso, ``Blockchain for transacting
  energy and carbon allowance in networked microgrids,'' \emph{IEEE Trans.
  Smart Grid}, vol.~12, no.~6, pp. 4702--4714, 2021.

\bibitem{zhong2023communication}
X.~Zhong, W.~Zhong, Y.~Liu, C.~Yang, and S.~Xie, ``A communication-efficient
  coalition graph game-based framework for electricity and carbon trading in
  networked energy hubs,'' \emph{Appl. Energy}, vol. 329, p. 120221, 2023.

\bibitem{yang2020distributed}
L.~Yang, J.~Luo, Y.~Xu, Z.~Zhang, and Z.~Dong, ``A distributed dual consensus
  admm based on partition for dc-dopf with carbon emission trading,''
  \emph{IEEE Trans. Ind. Informat.}, vol.~16, no.~3, pp. 1858--1872, 2020.

\bibitem{guo2021asynchronous}
Z.~Guo, P.~Pinson, Q.~Wu, S.~Chen, Q.~Yang, and Z.~Yang, ``An asynchronous
  online negotiation mechanism for real-time peer-to-peer electricity
  markets,'' \emph{IEEE Trans. Power Syst.}, pp. 1--13, 2021 early access.

\bibitem{2019zheng}
Y.~{Zheng}, Y.~{Song}, D.~J. {Hill}, and K.~{Meng}, ``Online distributed
  {MPC}-based optimal scheduling for {EV} charging stations in distribution
  systems,'' \emph{IEEE Trans. Ind. Informat.}, vol.~15, no.~2, pp. 638--649,
  Feb. 2019.

\bibitem{jiao2022online}
F.~Jiao, Y.~Zou, X.~Zhang, and B.~Zhang, ``Online optimal dispatch based on
  combined robust and stochastic model predictive control for a microgrid
  including {EV} charging station,'' \emph{Energy}, vol. 247, p. 123220, 2022.

\bibitem{fan2020online}
S.~Fan, G.~He, X.~Zhou, and M.~Cui, ``Online optimization for networked
  distributed energy resources with time-coupling constraints,'' \emph{IEEE
  Trans. Smart Grid}, vol.~12, no.~1, pp. 251--267, 2020.

\bibitem{2017RTMG}
W.~{Shi}, N.~{Li}, C.~{Chu}, and R.~{Gadh}, ``Real-time energy management in
  microgrids,'' \emph{IEEE Trans. Smart Grid}, vol.~8, no.~1, pp. 228--238,
  Jan. 2017.

\bibitem{zhong2019online}
W.~Zhong, K.~Xie, Y.~Liu, C.~Yang, S.~Xie, and Y.~Zhang, ``Online control and
  near-optimal algorithm for distributed energy storage sharing in smart
  grid,'' \emph{IEEE Trans. Smart Grid}, vol.~11, no.~3, pp. 2552--2562, 2020.

\bibitem{yu2018distributed}
L.~Yu, T.~Jiang, and Y.~Zou, ``Distributed real-time energy management in data
  center microgrids,'' \emph{IEEE Trans. Smart Grid}, vol.~9, no.~4, pp.
  3748--3762, 2018.

\bibitem{pjm}
{PJM-Data Miner 2}, ``Real-time five minute {LMPs},'' 2024, accessed Mar.,
  2024. [Online]. Available: https://dataminer2.pjm.com/.

\bibitem{mohamed2014real}
A.~Mohamed, V.~Salehi, T.~Ma, and O.~Mohammed, ``Real-time energy management
  algorithm for plug-in hybrid electric vehicle charging parks involving
  sustainable energy,'' \emph{IEEE Trans. Sustain. Energy}, vol.~5, no.~2, pp.
  577--586, 2014.

\bibitem{jin2014optimized}
C.~Jin, X.~Sheng, and P.~Ghosh, ``Optimized electric vehicle charging with
  intermittent renewable energy sources,'' \emph{IEEE J. Sel. Top. Signal
  Proc.}, vol.~8, no.~6, pp. 1063--1072, 2014.

\end{thebibliography}

\appendices
\makeatletter
\@addtoreset{equation}{section}
\@addtoreset{theorem}{section}
\makeatother
\setcounter{equation}{0}  
\renewcommand{\theequation}{A.\arabic{equation}}
\renewcommand{\thetheorem}{A.\arabic{theorem}}
\section{Proof of Proposition \ref{prop-1}}
\label{appendix-A}
Let $\{p_{d,t},\forall t\}$ be the aggregate power trajectory. For each time slot $t \in \mathcal{T}$, since $p_{d,t} \in [\check p_{d,t}^*, \hat p_{d,t}^*]$, we can define an auxiliary coefficient:
\begin{align}
    \beta_t:= \frac{\hat p_{d,t}^*-p_{d,t}}{\hat p_{d,t}^*-\check p_{d,t}^*} \in [0,1]
\end{align}
so that $p_{d,t}=\beta_t \check{p}^*_{d,t}+(1-\beta_t)\hat{p}^*_{d,t}$. Then, we can construct a feasible EV dispatch strategy by letting
\begin{subequations}
\begin{align}
    p^c_{v,t}&=\beta_t\check{p}^{c*}_{v,t}+(1-\beta_t)\hat{p}^{c*}_{v,t},\\
    e_{v,t}&=\beta_t\check{e}^{c*}_{v,t}+(1-\beta_t)\hat{e}^{c*}_{v,t}.
\end{align}
\end{subequations}
for all time slots $t \in \mathcal{T}$.

We prove that it is a feasible EV dispatch strategy as follows,
\begin{align}
    p_{d,t}=~ & \beta_t\check p_{d,t}^* + (1-\beta_t) \hat p_{d,t}^* \nonumber\\
    =~ &  \beta_t \sum_{v \in \mathcal{V}} \check p_{v,t}^{c*} + (1-\beta_t) \sum_{v \in \mathcal{V}} \hat p_{v,t}^{c*} \nonumber\\
    = ~ & \sum_{v \in \mathcal{V}} \left[\beta_t \check p_{v,t}^{c*}+(1-\beta_t) \hat p_{v,t}^{c*}\right] \nonumber\\
    = ~ & \sum_{v \in \mathcal{V}} p^c_{v,t}
\end{align}
Hence, constraint \eqref{equ:ubpdi} holds for $p_{d,t}$ and $p^c_{v,t},\forall v$. Similarly, we can prove that constraints \eqref{equ:ubpd}-\eqref{equ:uberanges} and \eqref{eq:evava} are met. Therefore, we have constructed a feasible EV dispatch strategy, which completes the proof. \hfill$\blacksquare$

\setcounter{equation}{0}  
\renewcommand{\theequation}{B.\arabic{equation}}
\renewcommand{\thetheorem}{B.\arabic{theorem}}
\section{Proof of Proposition \ref{prop-2}}
\label{appendix-B}
Denote the optimal solution of \textbf{P1'} as $\hat{x}_{g,t}^{*}$ and $\check{x}_{g,t}^{*}$, and the optimal solution of \textbf{P1} by $\hat{x}_{g,t}^{off}$ and $\check{x}_{g,t}^{off}$. According to \eqref{equ:dppInequ}, we have
\begin{equation}\label{eq:gap}
\begin{aligned}
&\mathbb{E}\left[\Delta(\boldsymbol{\Theta}_t)|\boldsymbol{\Theta}_t\right]+V_1\mathbb{E}[-F_t^{*}|\boldsymbol{\Theta}_t]\\
&\leq A_1+V_1\mathbb{E}[-F_t^{*}|\boldsymbol{\Theta}_t]+\sum\limits_{g\in{\mathcal{G}}}\hat{Q}_{g,t}\mathbb{E}\left[\hat{a}_{g,t}-\hat{x}_{g,t}^{*}|\boldsymbol{\Theta}_t\right]\\
&+\sum\limits_{g\in{\mathcal{G}}}\check{Q}_{g,t}\mathbb{E}\left[\check{a}_{g,t}-\check{x}_{g,t}^{*}|\boldsymbol{\Theta}_t\right]+\sum\limits_{g\in{\mathcal{G}}}\hat{Z}_{g,t}\mathbb{E}\left[-\hat{x}_{g,t}^{*}|\boldsymbol{\Theta}_t\right]\\
&+\sum\limits_{g\in{\mathcal{G}}}\check{Z}_{g,t}\mathbb{E}\left[-\check{x}_{g,t}^{*}|\boldsymbol{\Theta}_t\right],\\
&\leq A_1+V_1\mathbb{E}[-F_t^{off}|\boldsymbol{\Theta}_t]+\sum\limits_{g\in{\mathcal{G}}}\hat{Q}_{g,t}\mathbb{E}\left[\hat{a}_{g,t}-\hat{x}_{g,t}^{off}|\boldsymbol{\Theta}_t\right]\\
&+\sum\limits_{g\in{\mathcal{G}}}\check{Q}_{g,t}\mathbb{E}\left[\check{a}_{g,t}-\check{x}_{g,t}^{off}|\boldsymbol{\Theta}_t\right]+\sum\limits_{g\in{\mathcal{G}}}\hat{Z}_{g,t}\mathbb{E}\left[-\hat{x}_{g,t}^{off}|\boldsymbol{\Theta}_t\right]\\
&+\sum\limits_{g\in{\mathcal{G}}}\check{Z}_{g,t}\mathbb{E}\left[-\check{x}_{g,t}^{off}|\boldsymbol{\Theta}_t\right].
\end{aligned}
\end{equation}
By summing the above inequality \eqref{eq:gap} over time slots $t\in\{1,2,\ldots,T\}$, dividing both sides by $V_1 T$, and letting $T \to \infty$, we have
\begin{align}
   & \lim_{T \to \infty} \frac{1}{T}\left(\mathbb{E}[L(\boldsymbol{\Theta}_{T+1})]-\mathbb{E}[L(\boldsymbol{\Theta}_1)]\right) + \lim_{T \to \infty} \frac{V_1}{T} \sum_{t=1}^T \mathbb{E}(-F_t^{*}) \nonumber\\
  & \le  A_1 + \lim_{T \to \infty} \frac{V_1}{T} \sum_{t=1}^T \mathbb{E}(-F_t^{off}).
\end{align}

This is based on the fact that
\begin{align}
    \lim\limits_{T\rightarrow\infty}\frac{1}{T}\sum\limits_{t=1}^{T}\mathbb{E}\left[\hat{a}_{g,t}-\hat{x}_{g,t}^{off}|\boldsymbol{\Theta}_t\right]\leq 0,\\
    \lim\limits_{T\rightarrow\infty}\frac{1}{T}\sum\limits_{t=1}^{T}\mathbb{E}\left[\check{a}_{g,t}-\check{x}_{g,t}^{off}|\boldsymbol{\Theta}_t\right]\leq 0,\\
    \lim\limits_{T\rightarrow\infty}\frac{1}{T}\sum\limits_{t=1}^{T}\mathbb{E}\left[-\hat{x}_{g,t}^{off}|\boldsymbol{\Theta}_t\right]\leq 0,\\
    \lim\limits_{T\rightarrow\infty}\frac{1}{T}\sum\limits_{t=1}^{T}\mathbb{E}\left[-\check{x}_{g,t}^{off}|\boldsymbol{\Theta}_t\right]\leq 0,
\end{align}
which is due to constraints \eqref{eq:Qgub-lim}-\eqref{eq:xglb}.

Since $L(\boldsymbol{\Theta}_{T+1})$ and $L(\boldsymbol{\Theta}_1)$ are finite, we have
\begin{align*}
\underbrace{\lim\limits_{T\rightarrow\infty}\frac{1}{T}\sum\limits_{t=1}^{T}\mathbb{E}[F_t^{*}]}_{\mathcal{F}^{*}}\ge -\frac{A_1}{V_1}+\underbrace{\lim\limits_{T\rightarrow\infty}\frac{1}{T}\sum\limits_{t=1}^{T}\mathbb{E}[F_t^{off}]}_{\mathcal{F}^{off}}.
\end{align*}
Moreover, it is easy to know $\mathcal{F}^{off} \ge \mathcal{F}^{*}$. $\hfill \blacksquare$

\setcounter{equation}{0}  
\renewcommand{\theequation}{C.\arabic{equation}}
\renewcommand{\thetheorem}{C.\arabic{theorem}}
\section{Proof of Proposition \ref{prop-3}}
\label{appendix-C}
Suppose constraint \eqref{eq:c-lbub} holds in time slot $t$, we next prove that it also holds in time slot $t+1$ by induction.

\textbf{Case 1}: $0\leq c_{t} < m^{b,max}$.
The partial derivative of the objective function in $\textbf{P2}'$, denoted by $P2'_t$, with respect to $p^{g}_{t}$ is
\begin{align*}
\frac{\partial P2'_t}{\partial p^{g}_{t}}& = V_2\frac{\partial C_t}{\partial p^{g}_{t}} + H_{t}\rho_t\\
                                          &\leq V_2\pi^{g,max}_{t}+\left(c_{t} - \phi_{t}\right)\rho_t\\
                                          &=V_2\pi^{g,max}_{t}+\left(c_{t} - m^{b,max} - V_2\pi_t^{g,max}\frac{1}{\rho_t}\right)\rho_t\\
                                          &=\left(c_{t} - m^{b,max}\right)\rho_t<0.
\end{align*}
Thus, the objective function is strictly decreasing with respect to $p^{g}_{t}$.
Therefore, the optimal solution is $p^{g}_{t}=p^{g,max}$.
In addition,
\begin{align*}
\frac{\partial P2'_t}{\partial m_{t}^{b}} &= V_2\frac{\partial C_t}{\partial m^{b}_{t}} - H_{t}\\
&\geq V_2 \pi_t^{c,min} - (c_{t} - m^{b,max} - V_2\pi_t^{g,max}\frac{1}{\rho_t})\\
&=V_2(\pi_t^{c,min} + \pi_t^{g,max}\frac{1}{\rho_t}) + m^{b,max} - c_{t} > 0.
\end{align*}
the objective function is strictly increasing with respect to $m_t^b$. The optimal solution is $m^{b}_{t}=0$.
Further, based on \eqref{eq:c-state}, we have $c_{t+1} = c_{t}+\rho_t p^{g,max}$ and hence
\begin{align*}
0 \leq c_{t+1} \leq m^{b,max}+\rho_tp^{g,max} \leq c^{max}.
\end{align*}
The third inequality is due to the assumption in Proposition \ref{prop-3}.

\textbf{Case 2}: $m^{b,max} \leq c_{t} \leq V_2\left(\pi_t^{b,max}+\pi_t^{g,max}\frac{1}{\rho_t^{min}}\right) + m^{b,max}$.
Due to $V_2 \leq V_{2,max} \leq \frac{c^{max}_{t}-m^{b,max} - \rho_t^{max} p^{g,max}}{\pi_t^{b,max}+\pi_t^{g,max}\frac{1}{\rho^{min}}}$, we have
$c_{t} \leq c^{max}_{t} - \rho_t^{max} p^{g,max}$.
Thus, based on the update \eqref{eq:c-state}, we have
\begin{align*}
  c_{t+1} & \leq c^{max} - \rho_t^{max} p^{g,max} + \rho_t p^g_{t} - m^{b}_{t} \leq c^{max}.
\end{align*}
In addition, since $m^{b,max} \leq c_{t}$, we can obtain
\begin{align*}
  c_{t+1} & \geq m^{b,max} - m^{b}_{t} + \rho_t p^{g}_{t} \geq 0.
\end{align*}

\textbf{Case 3}: $V_2\left(\pi_t^{b,max}+\pi_t^{g,max}\frac{1}{\rho_t^{min}}\right) + m^{b,max} < c_{t} \leq c^{max}$.
Due to \eqref{eq:Vmax}, we have
$V_2\left(\pi_t^{b,max}+\pi_t^{g,max}\frac{1}{\rho_t^{min}}\right) + m^{b,max}
\leq c^{max}_{t} - \rho_t^{max} p^{g,max}
< c^{max}$.
Similar to \textbf{Case 1}, we then derive the partial derivative of the objective function of $P2'_t$ with respect to $m^{b}_{t}$, i.e.,
\begin{align*}
\frac{\partial P2'_t}{\partial m_{t}^{b}} &= V_2\frac{\partial C_t}{\partial m^{b}_{t}} - H_{t}\\
&\leq V_2 \pi_t^{c,max} - (c_{t} - m^{b,max} - V_2\pi_t^{g,max}\frac{1}{\rho_t})\\
&=V_2(\pi_t^{c,max} + \pi_t^{g,max}\frac{1}{\rho_t}) + m^{b,max} - c_{t}  <0.
\end{align*}
Thus, the objective function is strictly decreasing with respect to $m^{b}_{t}$.
Therefore, the optimal solution is $m^{b}_{t}=m^{b,max}$.
Meanwhile,
\begin{align*}
\frac{\partial P2'_t}{\partial p^{g}_{t}}& = V_2\frac{\partial C_t}{\partial p^{g}_{t}} + H_{t}\rho_t\\
                                        &\geq V_2\pi^{g,min}_{t} + \left(c_{t} - \phi_{t}\right)\rho_t\\
                                        & = V_2\pi^{g,min}_{t} + \left(c_{t} - m^{b,max} - V_2\pi_t^{g,max}\frac{1}{\rho_t}\right)\rho_t\\
                                        & = V_2\pi^{g,min}_{t} + V_2\pi_t^{b,max}\rho_t > 0.
\end{align*}
the objective function is strictly increasing with respect to $p^{g}_{t}$.
Therefore, the optimal solution is $p^{g}_{t}=0$.
According to \eqref{eq:c-state}, we have $c_{t+1}= c_{t} - m^{b,max}$ and hence
\begin{align*}
0 \leq c_{t+1} \leq c^{max} - m^{b,max} \leq c^{max}.
\end{align*}
Therefore, we have proved that the hard constraint \eqref{eq:c-state} still holds for all time slots. $\hfill\blacksquare$

\setcounter{equation}{0}  
\renewcommand{\theequation}{D.\arabic{equation}}
\renewcommand{\thetheorem}{D.\arabic{theorem}}
\section{Proof of Proposition \ref{prop-4}}
\label{appendix-D}
Denote $\widehat{m}_{c,t}$ and $\widehat{C}_t$ as the optimal results based on the optimal solution of $\textbf{P2}$ in time slot $t$.
Denote $m^*_{c,t}$ and $C_t^*$ as the optimal results of $\textbf{P2}'$ in time slot $t$.
According to \eqref{eq:dppInequ}, we have
\begin{align*}
&\Delta(H_t)+V_2\mathbb{E}[\widehat{C}_t|H_t]\\
& \leq A_2 + \sum\limits_{i\in{\mathcal{I}}}H_t\mathbb{E}\left[\widehat{m}_{c,t}|H_t\right] + V_2\mathbb{E}[\widehat{C}_t|H_t]\nonumber\\
& \leq A_2 + \sum\limits_{i\in{\mathcal{I}}}H_t\mathbb{E}\left[p^*_{c,t}|H_t\right] + V_2\mathbb{E}[C_t^*|H_t]\\
& = A_2 + \sum\limits_{i\in{\mathcal{I}}}H_t\mathbb{E}\left[m^*_{c,t}\right] + V_2\mathbb{E}[C_t^*].
\end{align*}
Since the system state is i.i.d., $m^*_{c,t}$ is also i.i.d. stochastic process.
Then, according to the strong law of large numbers, we obtain
\begin{align*}
&  \mathbb{E}[L(H_{t+1}) - L(H_t)|H_t]+V\mathbb{E}[\widehat{C}_t|H_t]\\
&  \leq A_2+\sum\limits_{i\in{\mathcal{I}}}H_t\lim\limits_{T\rightarrow\infty}\frac{1}{T}\sum\limits_{t=1}^{T}\left\{m^*_{c,t}\right\} + V\mathbb{E}[C_t^*].
\end{align*}
By taking the expectation of the above inequality, we have
\begin{align*}
&  \mathbb{E}[L(H_{t+1})] - \mathbb{E}[L(H_t)] + V_2\mathbb{E}(\widehat{C}_t)\\
&  \leq A_2 + \sum\limits_{i\in{\mathcal{I}}}H_t\lim\limits_{T\rightarrow\infty}\frac{1}{T}\sum\limits_{t=1}^{T}\mathbb{E}\left[m^*_{c,t}\right] + V_2\mathbb{E}[C_t^*].\\
& \leq A_2 + V_2\mathbb{E}[C_t^*].
\end{align*}
By summing the above inequality over time slots $t\in\{1,2,\ldots,T\}$, we have
\begin{equation*}
\sum\limits_{t=1}^{T}V_2\mathbb{E}[\widehat{C}_t] \leq A_2 T + V_2\sum\limits_{t=1}^{T}\mathbb{E}[C_t^*]-\mathbb{E}[L(H_{T+1})] + \mathbb{E}[L(H_1)].
\end{equation*}
Since $L(H_{T+1})$ and $L(H_{1})$ are finite, we divide both sides of the above inequalities by $V_2 T$ and take limits as $T\rightarrow\infty$ yielding
\begin{align*}
\lim\limits_{T\rightarrow\infty}\frac{1}{T}\sum\limits_{t=1}^{T}\mathbb{E}(\widehat{C}_t\leq \frac{A_2}{V_2}+\lim\limits_{T\rightarrow\infty}\frac{1}{T}\sum\limits_{t=1}^{T}\mathbb{E}(C_t^*).
\end{align*}
So far, we have finished the proof. $\hfill\blacksquare$

\end{document}